\documentclass{amsproc}

\usepackage{}

\newtheorem{theorem}{Theorem}[section]

\theoremstyle{definition}
\newtheorem{definition}[theorem]{Definition}
\newtheorem{example}[theorem]{Example}

\newtheorem{conjecture}[theorem]{Conjecture}

\newtheorem{proposition}[theorem]{Proposition}

\theoremstyle{remark}
\newtheorem{remark}[theorem]{Remark}

\numberwithin{equation}{section}

\DeclareMathOperator{\Aut}{Aut}

\DeclareMathOperator{\Spec}{Spec}

\DeclareMathOperator{\Def}{Def}

\DeclareMathOperator{\cHom}{\mathcal{H}\mathnormal{om}}
\DeclareMathOperator{\Ext}{Ext}
\DeclareMathOperator{\cExt}{\mathcal{E}\mathnormal{xt}}

\DeclareMathOperator{\Pic}{Pic}

\DeclareMathOperator{\PGL}{PGL}
\DeclareMathOperator{\SL}{SL}

\DeclareMathOperator{\GIT}{GIT}

\DeclareMathOperator{\End}{End}

\DeclareMathOperator{\rk}{rk}

\newcommand{\QED}{\ifhmode\unskip\nobreak\fi\quad {\rm Q.E.D.}} 

\newcommand{\bC}{\mathbb C}

\newcommand{\bP}{\mathbb P}
\newcommand{\bQ}{\mathbb Q}
\newcommand{\bR}{\mathbb R}
\newcommand{\bT}{\mathbb T}
\newcommand{\bZ}{\mathbb Z}

\newcommand{\cX}{\mathcal X}
\newcommand{\cY}{\mathcal Y}
\newcommand{\sD}{\mathcal D}
\newcommand{\cE}{\mathcal E}
\newcommand{\cF}{\mathcal F}

\newcommand{\sH}{\mathcal H}

\newcommand{\cO}{\mathcal O}

\newcommand{\sL}{\mathcal L}
\newcommand{\cM}{\mathcal M}
\newcommand{\cN}{\mathcal N}
\newcommand{\cP}{\mathcal P}

\newcommand{\cT}{\mathcal T}
\newcommand{\cU}{\mathcal U}

\newcommand{\ocM}{\overline{\mathcal M}}

\DeclareMathOperator{\cEnd}{\mathcal{E}\mathnormal{nd}}

\DeclareMathOperator{\spz}{sp}

\DeclareMathOperator{\QG}{QG}
\DeclareMathOperator{\PU}{PU}
\DeclareMathOperator{\U}{U}

\newcommand{\ocP}{\overline{\mathcal P}}

\begin{document}

\title{Compact moduli spaces of surfaces of general type}


\author{Paul Hacking}
\address{Department of Mathematics and Statistics,
Lederle Graduate Research Tower, Box 34515,
University of Massachusetts,
Amherst, MA 01003-9305, USA}
\curraddr{}
\email{hacking@math.umass.edu}
\thanks{The author was partially supported by NSF grant DMS-0968824. I would like to thank T.~Bridgeland, W.~Chen, A.~Corti, I.~Dolgachev, D.~Huybrechts, J.~Koll\'ar, J.~Tevelev, R.~Thomas, G.~Urz\'ua, and J.~Wahl for helpful discussions and correspondence.}


\date{\today}

\begin{abstract}
We give an introduction to the compactification of the moduli space of surfaces of general type introduced by Koll\'ar and Shepherd-Barron and generalized to the case of surfaces with a divisor by Alexeev. The construction is an application of Mori's minimal model program for $3$-folds. We review the example of the projective plane with a curve of degree $d \ge 4$.
We explain a connection between the geometry of the boundary of the compactification of the moduli space and the classification of vector bundles on the surface in the case $H^{2,0}=H^1=0$.
\end{abstract}

\maketitle

\section{Introduction}

The moduli space $\cM_g$ of curves of genus $g \ge 2$ admits a compactification $\ocM_g$ given by the moduli space of stable curves \cite{DM69}.
Here a \emph{stable curve} $C$ is a connected compact curve of arithmetic genus $g$ having only nodal singularities $(xy=0) \subset \bC^2$ and finite automorphism group.
Components of the boundary $\partial\ocM_g := \ocM_g \setminus \cM_g$ correspond to the topological types of degenerations of a curve of genus $g$ obtained by contracting a simple loop to a point. The space $\ocM_g$ is an orbifold of dimension $3g-3$ and the boundary is a normal crossing divisor, that is, the pair $(\ocM_g,\partial\ocM_g)$ is locally isomorphic to $\bC^{3g-3}$ together with a collection of coordinate hyperplanes modulo a finite group.

The moduli space $\cM_{K^2,\chi}$ of surfaces of general type with fixed topological invariants admits an analogous compactification $\ocM_{K^2,\chi}$, the moduli space of \emph{stable surfaces} \cite{KSB88}. The definition of a stable surface is derived from Mori's minimal model program for $3$-folds.
Indeed, given a family of smooth surfaces of general type over a punctured disc,  results of the minimal model program produce (after a finite base change) a distinguished extension of the family over the disc (the \emph{canonical model} given by the projective spectrum of the canonical ring of a smooth model of the total space with special fiber a reduced normal crossing divisor). This is the exact analogue of the stable reduction theorem for curves, which uses the classical theory of minimal models of surfaces, cf. \cite[1.12]{DM69}. A \emph{stable surface} is then, roughly speaking, a surface which can arise as a limit of smooth surfaces in this way. The moduli space of stable surfaces gives a compactification of the moduli space of smooth surfaces, because it satisfies the valuative criterion of properness by construction.

Very little is known about the moduli spaces $\cM_{K^2,\chi}$ in general.
In particular they can be essentially arbitrarily singular \cite{V06}, and have arbitrarily many connected components \cite{C86}.
However, we expect that there are many examples of connected components of the moduli space that are well behaved.
A particularly attractive example is the case of surfaces of general type such that $H^{2,0}=0$ and $\pi_1=0$ discussed in \S\ref{conjecture}.

We want to study the compactification given by moduli of stable surfaces in these cases. Perhaps the most basic question is: what are the codimension $1$ components of the boundary? Here by definition the \emph{boundary} is the complement of the locus of surfaces with Du Val singularities (as these are the canonical models of smooth surfaces of general type). An important type of boundary divisor is the locus of equisingular deformations of a normal surface with a unique cyclic quotient singularity of type $\bC^2_{u,v}/(\bZ/n^2\bZ)$, where the action is given by
$$(u,v) \mapsto (\zeta u, \zeta^{na-1}v), \quad \zeta=\exp(2\pi i/n^2)$$
for some $n,a$ with $(a,n)=1$.
These singularities were first studied by J.~Wahl \cite[5.9.1]{W81} so we refer to them as \emph{Wahl singularities}.

The smoothing component of the deformation space of a Wahl singularity relevant to moduli of stable surfaces is a smooth curve germ.
So, at least in the absence of local-to-global obstructions, Wahl singularities occur in codimension $1$.
The Milnor fiber of the smoothing is a rational homology ball. Thus if $Y \leadsto X$ is a global degeneration which is locally of this form, specialization defines an isomorphism of rational homology $H_*(Y,\bQ) \rightarrow H_*(X,\bQ)$.
For this reason it seems that it is not possible to predict the existence of the degeneration $Y \leadsto X$ based on the topology of $Y$ alone.

In this paper we explain a relation between boundary components of the moduli space of Wahl type and the classification of rigid holomorphic vector bundles on the smooth surface $Y$ in the case $H^{2,0}=H^1=0$ \cite{H11}. We also describe in some detail a motivating example: the case of the projective plane $Y=\bP^2$ together with a curve $D$ of degree $d \ge 4$ \cite{H04}. This is an instance of the moduli space of stable surfaces with boundary, analogous to the moduli spaces of stable pointed curves $\ocM_{g,n}$. In this example we are able to describe the components of the boundary of the moduli space fairly explicitly. In particular, the boundary divisors associated to Wahl singularities display a beautiful combinatorial structure that is intimately related to the classification of rigid vector bundles on the projective plane.

\medskip
\noindent
\emph{Notation}: We work over the complex numbers.
We write $\bC^r/\frac{1}{n}(a_1,\ldots,a_r)$ or just $\frac{1}{n}(a_1,\ldots,a_r)$ for the cyclic quotient singularity
$\bC^r/(\bZ/n\bZ)$
where the action is given by
$$(x_1,\ldots,x_r) \mapsto (\zeta^{a_1}x_1,\ldots,\zeta^{a_r}x_r), \quad \zeta=\exp(2\pi i/n).$$

\section{Moduli of stable surfaces}

\subsection{Semi log canonical singularities}

Let $X$ be a reduced Cohen-Macaulay surface.
We assume that the set $S \subset X$ of singularities which are not of \emph{double normal crossing} type $(xy=0) \subset \bC^3$ is finite.
The \emph{double curve} $\Delta \subset X$ of $X$ is the union of the one dimensional components of the singular locus of $X$.
Let $\nu \colon X^{\nu} \rightarrow X$ be the normalization of $X$ and $\Delta^{\nu} \subset X^{\nu}$ the inverse image of the double curve. Thus the map $\Delta^{\nu} \rightarrow \Delta$ is  finite of degree $2$.
The dualizing sheaf $\omega_X$ of the surface $X$ is the correct generalization of the canonical line bundle for smooth surfaces.
It has the following properties: (1) $\omega_X$ is invertible away from the finite set $S$, (2) $\omega_X$ is reflexive, (3) $(\nu^*\omega_X)^{\vee\vee}=\omega_{X^{\nu}}(\Delta^{\nu})$. Here $\cF^{\vee\vee}$ denotes the double dual (or \emph{reflexive hull}) of a coherent sheaf $\cF$.

The dualizing sheaf $\omega_X$ corresponds to a Weil divisor class denoted $K_X$. Of course if $X$ is normal this is the usual canonical divisor class.

Let $X$ be a normal surface and $D$ a $\bQ$-divisor on $X$.
We say the pair $(X,D)$ is \emph{log canonical} if $K_X+D$ is $\bQ$-Cartier and the following condition holds.
Let $\pi \colon \tilde{X} \rightarrow X$ be a resolution of singularities of $X$ such that the sum of the exceptional divisor $E= \sum E_i$ of $\pi$ and the strict transform $D'$ of $D$ is a normal crossing divisor. Then
$$K_{\tilde{X}}+D' = \pi^*(K_{X}+D)+\sum a_iE_i$$
in $\Pic(\tilde{X}) \otimes \bQ$ for some $a_i \in \bQ$. We require that $a_i \ge -1$ for each $i$. (It suffices to check the condition for one such resolution.) See \cite[4.1]{KM98} for a detailed classification of log canonical surface-divisor pairs.

\begin{definition}\label{def-slc}
Let $X$ be a surface. We say $X$ has \emph{semi log canonical} (slc) singularities if the following conditions hold.
\begin{enumerate}
\item The surface $X$ is reduced, Cohen-Macaulay, and has only double normal crossing singularities $(xy=0) \subset \bC^3$ away from a finite set of points.
\item With notation as above, the pair $(X^{\nu},\Delta^{\nu})$ consisting of the normalization of $X$ together with the inverse image of the double curve has log canonical singularities.
\item For some $N>0$ the $N$th reflexive tensor power $\omega_X^{[N]}:= (\omega_X^{\otimes N})^{\vee\vee}$ of the dualizing sheaf $\omega_X$ is invertible.
\end{enumerate}
\end{definition}

Semi log canonical surface singularities were classified in \cite[\S4]{KSB88}. We describe the two examples we need in \S\ref{wahl} and \S\ref{nc} below.

\begin{remark}
The definition of semi log canonical singularities for $X$ of dimension $d$ is the same except that the Cohen-Macaulay condition is replaced by Serre's condition $S_2$ and we require that $X$ has only double normal crossing singularities $(xy=0) \subset \bC^{d+1}$ away from a codimension $2$ subset.
\end{remark}

\subsection{Stable surfaces}

\begin{definition} A \emph{stable surface} is a connected projective surface $X$ such that $X$ has semi log canonical singularities and the dualizing sheaf $\omega_X$ is ample.
\end{definition}

\subsection{Index one cover}

Let $(P \in X)$ be an slc surface germ.
The \emph{index} of $P \in X$ is the least integer such that $\omega_X^{[N]}$ is invertible.
Fix an isomorphism $\theta \colon \omega_X^{[N]} \rightarrow \cO_X$ and define a morphism $p \colon Z \rightarrow X$ by
$$Z:=\underline{\Spec}_X(\cO_X \oplus \omega_X \oplus \omega_X^{[2]} \oplus \cdots \oplus \omega_X^{[N-1]})$$
where the multiplication on $\cO_Z$ is defined using $\theta$.
Then $p \colon Z \rightarrow X$ is a cyclic cover of degree $N$ such that (1) the inverse image of $P$ is a single point $Q \in Z$, (2) the morphism $p$ is \'etale over $X \setminus \{P\}$, and (3) the surface $Z$ is Gorenstein, that is, $Z$ is Cohen-Macaulay and the dualizing sheaf $\omega_Z$ is invertible. Moreover, the germ $(Q \in Z)$ is slc.
The covering $p$ is called the \emph{index one cover} of $(P \in X)$. It is uniquely determined locally in the analytic topology.




\subsection{$\bQ$-Gorenstein deformations}\label{QG}

Let $(P \in X)$ be a slc surface germ. We say a deformation $(P \in \cX)/(0 \in S)$ is \emph{$\bQ$-Gorenstein} if it is induced by an equivariant deformation of the index one cover of $(P \in X)$.

Equivalently, let $\omega_{X/S}$ denote the relative dualizing sheaf of $X/S$. The sheaf $\omega_{X/S}$ is flat over $S$ and commutes with base change because $X/S$ is flat with Cohen-Macaulay fibers \cite[3.6.1]{C00}.
Then $X/S$ is $\bQ$-Gorenstein iff every reflexive power $\omega_{X/S}^{[N]}:=(\omega_{X/S}^{\otimes N})^{\vee\vee}$ of the relative dualizing sheaf $\omega_{X/S}$ commutes with basechange \cite[3.5]{H04}. (The definition of $\bQ$-Gorenstein deformations first appeared in this form in \cite[\S 6]{K91}.)

\subsection{Moduli stack}

\begin{definition}
Let $S$ be a scheme of finite type over $\bC$. A \emph{family of stable surfaces} over $S$ is a flat family $X/S$ such that each fiber is a stable surface and $X/S$ is $\bQ$-Gorenstein in the sense of \S\ref{QG}, that is, everywhere locally on $X$ the family $X/S$ is induced by an equivariant deformation of the index one cover of the fiber. Let $\ocM$ denote the moduli stack of stable surfaces so defined. For $K^2,\chi \in \bZ$ let $\ocM_{K^2,\chi}$ denote the moduli stack of stable surfaces $X$ with $K_X^2=K^2$, $\chi(\cO_X)=\chi$ (thus $\ocM_{K^2,\chi}$ is a union of connected components of $\ocM$).
\end{definition}

\begin{theorem}\cite{KSB88},\cite{K90},\cite{AM04}.
The stack $\ocM_{K^2,\chi}$ is a proper Deligne-Mumford stack of finite type over $\bC$ with projective coarse moduli space.
\end{theorem}

\section{Infinitesimal study of the moduli stack}

Let $X$ be a stable surface.
By definition, the moduli stack $\ocM$ of stable surfaces near $[X]$ is identified with the quotient stack $[\Def^{\QG}(X)/\Aut X]$
of the versal $\bQ$-Gorenstein deformation space $\Def^{\QG}(X)$ of $X$ by the automorphism group of $X$ (a finite group).
In this section we explain how to compute $\Def^{\QG}(X)$.

Let $X$ be a slc surface.
Recall that the $\bQ$-Gorenstein deformations of $X$ are by definition those deformations which are locally induced by an equivariant deformation of the index one cover. To treat global $\bQ$-Gorenstein deformations in a robust way it is best to rephrase this definition as follows.

The data of index one covers everywhere locally on $X$ defines a Deligne-Mumford stack $\mathfrak{X}$ with coarse moduli space $X$, the \emph{index one covering stack}.
Let $q \colon \mathfrak{X} \rightarrow X$ denote the induced map.
If $p \colon V \rightarrow U$ is an index one cover of an open set $U \subset X$ with group $G \simeq \bZ/N\bZ$, then
$q \colon \mathfrak{X} \rightarrow X$ is locally given by the map $[V/G] \rightarrow U=V/G$ from the quotient stack $[V/G]$ to its coarse moduli space. The deformations of $\mathfrak{X}$ are identified with the $\bQ$-Gorenstein deformations of $X$ \cite[3.7]{H04}, \cite{AH09}.
Note that the map $q \colon \mathfrak{X} \rightarrow X$ is an isomorphism away from the finite set of points $S \subset X$ where the dualizing sheaf $\omega_X$ is not invertible, that is, the stack structure is trivial away from $S$.

We consider coherent sheaves $\cF$ on the \'etale site of $\mathfrak{X}$. If $p$ is a local index one cover as above, the restriction of $\cF$ to $V$ is a $G$-equivariant sheaf $\cF_V$ on $V$, and $\cF$ is determined by these restrictions together with gluing data on the overlaps.
The pushforward $q_*\cF$ of $\cF$ may be defined locally by $(q_*\cF)|_U=(p_*\cF_V)^{G}$. Note that $q_*$ is exact.

The deformations of $\mathfrak{X}$ may be described using the cotangent complex $L_{\mathfrak{X}/\bC}$ for the structure sheaf $\cO_{\mathfrak{X}}$ on the \'etale site of $\mathfrak{X}$ \cite{I71}, \cite{I72}. We refer to \cite[\S3]{H10} for an exposition of the local theory of the cotangent complex following \cite{LS67}.
We define sheaves $\cT^i_{\QG,X}:=q_*\cExt^i(L_{\mathfrak{X}/\bC},\cO_{\mathfrak{X}})$ on $X$ and $\bC$-vector spaces $\bT^i_{\QG,X}:=\Ext^i(L_{\mathfrak{X}/\bC},\cO_{\mathfrak{X}})$. Then the basic result is the following.
\begin{proposition}
Let $X$ be a slc surface.
\begin{enumerate}
\item The first order infinitesimal $\bQ$-Gorenstein deformations of $X$ are identified with the $\bC$-vector space $\bT^1_{\QG,X}$.
\item Let $A$ be a finitely generated Artinian local $\bC$-algebra. Let $\cX/A$ be a $\bQ$-Gorenstein deformation of $X$ over $A$. Let $A' \rightarrow A$ be an infinitesimal extension of $A$ with kernel $I$. Let $\mathfrak{m}'$ denote the maximal ideal of $A'$ and assume that $\mathfrak{m}' \cdot I=0$ (so $I$ is a $A/\mathfrak{m}'=\bC$-vector space). Then there is an obstruction class $o(\cX/A,A') \in \bT^2_{\QG,X} \otimes I$ such that $o(\cX/A,A')=0$ iff there exists a $\bQ$-Gorenstein deformation $\cX'/A'$ of $\cX/A$ over $A'$. If $o(\cX/A,A')=0$ then isomorphism classes of such deformations form a torsor for $\bT^1_{\QG,X} \otimes I$.
\end{enumerate}
\end{proposition}
In particular, if $X$ is compact, the versal $\bQ$-Gorenstein deformation space of $X$ is an analytic germ $(0 \in \Def^{\QG}(X))$ with tangent space $\bT^1_{\QG,X}$ and is smooth if $\bT^2_{\QG,X}=0$.

Let $p \colon V \rightarrow U$ be an index one cover of an open set $U \subset X$ with group $G$.
Then we have $\cT^i_{\QG,X}|_U=(p_*\cT^i_V)^{G}$ where $\cT^i_V=\cExt^i(L_{V/\bC},\cO_V)$. We recall that, for any reduced scheme $V$, $\cT^0_V=\cT_V:=\cHom(\Omega_V,\cO_V)$ is the tangent sheaf of $V$, $\cT^1_V=\cExt^1(\Omega_V,\cO_V)$ and $\cT^2_V=0$ if $V$ is a local complete intersection. See for example \cite[p.~33]{A76}, \cite[3.6, 4.13]{H10}. It follows that
(1) $\cT^0_{\QG,X}=\cT_X$ is the tangent sheaf of $X$, (2) $\cT^1_{\QG,X}$ is supported on the singular locus of $X$ and can be computed locally by $\cT^1_{\QG,X}|_U=(p_*\cExt^1(\Omega_V,\cO_V))^{G}$, (3) $\cT^2_{\QG,X}$ is supported on the locus where the index one cover is not a local complete intersection.

We have the local-to-global spectral sequence
$$E^{p,q}_2=H^p(\cT^q_{\QG,X}) \Rightarrow \bT^{p+q}_{\QG,X}.$$
In particular, we obtain the exact sequence
$$0 \rightarrow H^1(\cT_X) \rightarrow \bT^1_{\QG,X} \rightarrow H^0(\cT^1_{\QG,X}) \rightarrow H^2(\cT_X).$$
Moreover, the obstruction space $\bT^2_{\QG,X}$ vanishes if $H^2(\cT_X)=0$, $H^1(\cT^1_{\QG,X})=0$, and the local index one covers are complete intersections.

\section{Wahl singularities}\label{wahl}

Let $P \in X$ denote the cyclic quotient singularity $\bC^2_{u,v}/\frac{1}{n^2}(1,na-1)$, where $(a,n)=1$.
These singularities were first studied by J.~Wahl \cite[5.9.1]{W81} so we refer to them as \emph{Wahl singularities}.
The index one cover $(Q \in Z) \rightarrow (P \in X)$ is the Du Val singularity of type $A_{n-1}$
$$\bC^2_{u,v} / \textstyle{\frac{1}{n}(1,-1)}=(xy=z^n) \subset \bC^3_{x,y,z}$$
where $x=u^n,y=v^n,z=uv$. Thus
$$(P \in X) =(xy=z^n) \subset \textstyle{\frac{1}{n}}(1,-1,a)$$
and a $\bQ$-Gorenstein smoothing is given by
\begin{equation}\label{Wahlsmoothing}
(P \in \cX)=(xy=z^n+t) \subset \textstyle{\frac{1}{n}}(1,-1,a) \times \bC^1_t.
\end{equation}

The one parameter deformation $(P \in \cX)/(0 \in \bC^1_t)$ is the versal $\bQ$-Gorenstein deformation of $(P \in X)$.
In particular, if $X$ is a stable surface with a unique singularity $P \in X$ of Wahl type and there are no local-to-global obstructions to $\bQ$-Gorenstein deformations of $X$, then the moduli stack of stable surfaces $\ocM$ is smooth near $[X]$ and equisingular deformations of $X$ determine a codimension $1$ component of the boundary.

\subsection{Topology of the smoothing}
If $P \in X$ is a normal surface singularity then the \emph{link} $L$ of $P \in X$ is the smooth $3$-manifold obtained as the intersection of $X$ with a small sphere about the point $P$ in some embedding $X \subset \bC^N$. For example the link of the cyclic quotient singularity
$\frac{1}{r}(1,a)$ is the lens space $L=L(r,a)=S^3/\frac{1}{r}(1,a) \subset \bC^2/\frac{1}{r}(1,a)$.

If $(P \in \cX)/(0 \in T)$ is a one-parameter smoothing of a normal surface singularity $P \in X$,
the \emph{Milnor fiber} $M$ of the smoothing is defined as follows. Choose an embedding $\cX \subset \bC^N \times T$, say $P \mapsto (0,0)$, and
let $B(\delta):=\{\|z\| \le \delta \} \subset \bC^N$ be the closed ball with center $0 \in \bC^N$ and radius $\delta$. We define
$M=\cX_{\epsilon} \cap B(\delta)$ for $0 < \epsilon \ll \delta \ll 1$. Then $M$ is a smooth $4$-manifold with boundary the link $L$ of the singularity $P \in X$.

Wahl observed that the Milnor fiber $M$ of the smoothing (\ref{Wahlsmoothing}) is a rational homology ball. More precisely, $\pi_1(M)$ is cyclic of order $n$ and $H_i(M,\bZ)=0$ for $i>1$. Indeed, the Milnor fiber $M$ is the quotient of the Milnor fiber $M_Z$ of the smoothing of the $A_{n-1}$ singularity $(Q \in Z)$ by a free $\bZ/n\bZ$-action. Now $M_Z$ is homotopy equivalent to a bouquet of $n-1$ copies of the $2$-sphere, in particular $M_Z$ is simply connected and $e(M_Z)=n$. Hence $\pi_1(M)=\bZ/n\bZ$ and $e(M)=e(M_Z)/n=1$. Finally $M$ is a Stein manifold of complex dimension $2$ so has the homotopy type of a cell complex of real dimension $2$. We deduce that $H_i(M)=0$ for $i>1$ as claimed.

The boundary of $M$ is the link $L=S^3/\frac{1}{n^2}(1,na-1)$ of $P \in X$ and the map $\pi_1(L) \rightarrow \pi_1(M)$ is the surjection $\bZ/n^2\bZ \rightarrow \bZ/n\bZ$.

Let $Y$ be a smooth compact surface. We would like to determine the possible degenerations of $Y$ to a singular surface $X$ which are locally of the form (\ref{Wahlsmoothing}). It seems that topological considerations provide very little information because the topology of $M$ is almost trivial.
We return to this question in \S\S\ref{bundles}--\ref{donaldson} below, where we describe an approach in the case $H^{2,0}(Y)=H^1(Y)=0$.

\section{Orbifold double normal crossing singularities}\label{nc}

Let $X$ be a surface which is a union of two smooth components $X_1$, $X_2$ meeting transversely along a curve $C$.
That is, $X$ has double normal crossing singularities $(xy=0) \subset \bC^3$ along $C$.
Then the sheaf $\cT^1_X$ is the line bundle on $C$ given by the tensor product of the normal bundles of $C$ in $X_1$ and $X_2$
$$\cT^1_X=\cExt^1(\Omega_X,\cO_X)= \cN_{C/X_1} \otimes \cN_{C/X_2}.$$
See \cite[2.3]{F83}. If $U \subset X$ is a local chart of the form $(xy=0) \subset \bC^3$ then $\cT^1_X|_U$ is generated by $x^{-1} \otimes y^{-1}$, and a local section $s=f(z) \cdot x^{-1} \otimes y^{-1} \in H^0(U,\cT^1_X)$ corresponds to the first order infinitesimal deformation
$$\cU=(xy= t \cdot f(z)) \subset \bC^3_{x,y,z} \times (\Spec \bC[t]/(t^2)).$$
A necessary condition for smoothability of $X$ is $H^0(\cT^1_X) \neq 0$ (because otherwise all deformations are locally trivial).
Suppose $\cX/(0 \in \bC^1_t)$ is a one-parameter deformation of $X$ and $s \in H^0(\cT^1_X)$ is the induced element (the local part of the induced first order infinitesimal deformation). Then the total space $\cX$ is smooth away from the zeroes of $s$ and has ordinary double points $(xy+zt=0)$ at simple zeroes of $s$. In particular, if $s$ has only simple zeroes then the general fiber is smooth.

The \emph{orbifold double normal crossing singularity}
$$(xy=0) \subset \bC^3_{x,y,z} / \textstyle{\frac{1}{n}}(1,-1,a), \quad (a,n)=1,$$
is semi log canonical and occurs frequently on stable limits of smooth surfaces.
The index one cover is the double normal crossing singularity $\mbox{$(xy=0)$} \subset \bC^3$. A $\bQ$-Gorenstein smoothing is given by
$$(xy=t) \subset \textstyle{\frac{1}{n}}(1,-1,a) \times \bC^1_t.$$

Let $X$ be a surface that is a union of two normal components $X_1$ and $X_2$ meeting along a smooth curve $C$ such that $X$ has orbifold double normal crossing singularities along $C$ and is smooth elsewhere. Then the $\bQ$-Gorenstein deformations of $X$ are, by definition, locally induced by equivariant deformations of the index one cover, which is a double normal crossing singularity. The sheaf $\cT^1_{\QG,X}$ of local first order $\bQ$-Gorenstein deformations is a line bundle on $C$ and may be computed as follows.
Let $C_i|_C$ denote the $\bQ$-divisor on $C$, well defined up to linear equivalence, obtained by moving $C$ on $X_i$ and restricting to $C$.
Then the sum $C_1|_C+C_2|_C$ is a $\bZ$-divisor on $C$, and
\begin{equation}\label{QGinfnorbdle}
\cT^1_{\QG,X}=\cO_C(C_1|_C+C_2|_C).
\end{equation}
In particular $\cT^1_{\QG,X}$ is a line bundle on $C$ of degree $(C^2)_{X_1}+(C^2)_{X_2}$ equal to the sum of the self-intersection numbers of $C$ on $X_1$ and $X_2$.
If $U \subset X$ is a local chart of the form $(xy=0) \subset \frac{1}{n}(1,-1,a)$, then $\cT^1_{\QG,X}|_U$ is generated by $x^{-1} \otimes y^{-1}$,
and a local section $s=f(z^n) \cdot x^{-1} \otimes y^{-1}$ corresponds to the first order deformation
$$(xy=t \cdot f(z^n)) \subset \textstyle{\frac{1}{n}}(1,-1,a) \times \Spec(\bC[t]/(t^2)).$$
In particular, if $s$ has a simple zero at the orbifold point $(P \in C) = (0 \in \bC^1_{z^n})$ then locally this deformation may be extended to a one-parameter deformation with general fiber the Wahl singularity $\frac{1}{n^2}(1,na-1)$.

\section{Surfaces with boundary}\label{modulipairs}

The moduli space of pairs $(X,D)$ consisting of a smooth surface $X$ and an effective $\bQ$-divisor $D$ such that $(X,D)$ is log canonical and $K_X+D$ is ample admits a compactification given by the moduli space of stable pairs \cite{A96}, \cite{K11}. These compactifications are analogous to the moduli spaces of pointed stable curves $\ocM_{g,n}$ and their generalizations given by assigning weights to the marked points \cite{H03}.

\begin{definition}
Let $(X,D)$ be a pair consisting of a surface $X$ and an effective $\bQ$-divisor $D$.
We say the pair $(X,D)$ is a \emph{stable pair} or \emph{stable surface with boundary} if the following conditions hold.
\begin{enumerate}
\item The surface $X$ is reduced, Cohen-Macaulay, and has double normal crossing singularities away from a finite set. No component of the double curve is contained in the support of the divisor $D$.
\item The pair $(X^{\nu},\Delta^{\nu}+D^{\nu})$, consisting of the normalization of $X$ together with the sum of the inverse images of the double curve and the divisor $D$, is log canonical.
\item The $\bQ$-divisor $K_X+D$ is $\bQ$-Cartier and ample.
\end{enumerate}
\end{definition}


The definition of a family of stable pairs $(X,D)$ over a scheme $S$ involves technical difficulties in general \cite{A08}, \cite[\S6]{K10}.
These difficulties do not occur if $D$ is a $\bZ$-divisor or if we use \emph{floating coefficients} $D=(\alpha+\epsilon)B$ where $\alpha \in [0,1) \cap \bQ$ and $0 < \epsilon \ll 1$, $\epsilon \in \bQ$. We give a precise definition in the example described in \S\ref{planecurves} below, where floating coefficients are used.


\section{Plane curves}\label{planecurves}

A compactification of the moduli space of plane curves was obtained in the author's PhD thesis advised by A.~Corti \cite{H04}.
It is an instance of the moduli space of stable surfaces with boundary.
In this section we review the description of this example as it displays many features we expect to observe in the general case.

Let $\cP_d$ denote the moduli stack of smooth plane curves of degree $d \ge 4$.
Explicitly $\cP_d$ is the quotient stack $[U / \PGL(3)]$ where $U \subset \bP^N$ is the complement of the discriminant locus in the projective space of homogeneous polynomials of degree $d$ in $3$ variables, $N=\frac{1}{2}d(d+3)$.
We regard $\cP_d$ as the moduli stack of pairs $(X,D)$ consisting of a surface $X$ isomorphic to $\bP^2$ and a divisor $D$ on $X$ given by a smooth plane curve of degree $d$.
The pair $(X,\alpha D)$ has log canonical singularities and $K_X+\alpha D$ is ample for $3/d< \alpha \le 1$, $\alpha \in \bQ$.
One can thus define a compactification $\ocP^{\alpha}_d$ of $\cP_d$ by taking the closure of the image of $\cP_d$ in the moduli stack of stable pairs under the map $[D] \mapsto [(X,\alpha D)]$. Roughly speaking, the compactification $\ocP^{\alpha}_d$ becomes more complicated as $\alpha$ increases. So our approach is to consider the case where $\alpha$ is ``as small as possible". More precisely, there exists $\epsilon_0 > 0$ such that our compactification $\ocP_d$ coincides with $\ocP^{\alpha}_d$ for $3/d < \alpha < 3/d +\epsilon_0$.

As noted in \S\ref{modulipairs}, the definition of the moduli stack of stable pairs in general is rather technical. So we define the moduli stack $\ocP_d$ directly as follows.

\begin{definition}
A \emph{stable plane pair} of degree $d$ is a pair $(X,D)$ consisting of a surface $X$ together with an effective Weil divisor $D$ such that the pair $(X,(3/d+\epsilon)D)$ has semi log canonical singularities and the $\bQ$-Cartier divisor $K_X+(3/d+\epsilon)D$ is ample for $0< \epsilon \ll 1$, $\epsilon \in \bQ$.
\end{definition}

Let $(X,D)$ be a stable plane pair. Then in particular $X$ is slc and $D$ is $\bQ$-Cartier.
Let $P \in X$ be a point and $p \colon (Q \in Z) \rightarrow (P \in X)$ be the index one cover.
Let $D_Z$ be the divisor on $Z$ given by the inverse image of $D$ on $X$. We say $(X,D)$ is \emph{pre-smoothable} if $D_Z$ is Cartier for each $P \in X$. (This condition simplifies the deformation theory of pairs considered below. One can show that it is satisfied by the stable limits of pairs $(\bP^2,C)$ consisting of the plane together with a curve \cite[3.13]{H04}.)

Let $(X,D)$ be a pre-smoothable stable plane pair.
We regard $D$ as a closed subscheme $D \subset X$ of codimension $1$ (without embedded points).
By a deformation of the pair $(X,D)$ over a germ $(0 \in S)$ we mean a pair $(\cX,\sD)/(0 \in S)$ consisting of a deformation $\cX/(0 \in S)$ of $X$ together with a closed subscheme $\sD \subset \cX$ such that $\sD_0=D$ and $\sD$ is flat over $S$.
We say a deformation of $(P \in X,D)$ is \emph{$\bQ$-Gorenstein} if it is induced by an equivariant deformation of $(Q \in Z, D_Z)$.

A \emph{family of pre-smoothable stable plane pairs} over a scheme $S$ is a flat family $(\cX,\sD)/S$ such that each fiber is a pre-smoothable stable plane pair and the family is $\bQ$-Gorenstein, that is, everywhere locally on $\cX$ induced by an equivariant deformation of the index one cover of the fiber.
Let $\tilde{\cP}$ denote the stack of pre-smoothable stable plane pairs so defined. Then $\tilde{\cP}$ is a Deligne-Mumford stack locally of finite type over $\bC$. (Note: we do \emph{not} assert that $\tilde{\cP}$ is proper or that its connected components are of finite type.)
Let $\ocP_d \subset \tilde{\cP}$ denote the closure of $\cP_d$, the stack of \emph{smoothable stable plane pairs} of degree $d$.

\begin{theorem}
The stack $\ocP_d$ is a proper Deligne-Mumford stack of finite type over $\bC$. If $d$ is not a multiple of $3$ then $\ocP_d$ is smooth and is a connected component of $\tilde{\cP}$.
\end{theorem}

\subsection{Classification of stable planes}

\begin{proposition}
Let $(X,D)$ be a smoothable stable plane pair of degree $d$. Then $-K_X$ is ample and $dK_X+3D$ is linearly equivalent to zero.
\end{proposition}

We sketch the proof of the proposition. Since $[(X,D)]$ lies in the closure of $\cP_d$, there is a $\bQ$-Gorenstein deformation $(\cX,\sD)/(0 \in T)$ of $(X,D)$ over the germ of a curve such that the general fiber is the plane together with a curve of degree $d$.
The restriction of $dK_{\cX}+3\sD$ to a general fiber of the family $\cX/(0 \in T)$ is linearly equivalent to zero,
hence $dK_{\cX}+3\sD$ is linearly equivalent to a sum of components of the special fiber $X$. But also $dK_X+3D=(dK_{\cX}+3\sD)|_X$ is nef because $K_{X}+(3/d+\epsilon)D$ is ample for $0<\epsilon \ll 1$. It follows that $dK_X+3D \sim 0$.
Now $-K_X$ is ample because $K_X+(3/d+\epsilon)D$ is ample.

If $d$ is not divisible by $3$ then, since $dK_X+3D \sim 0$, the canonical class $K_X$ is divisible by $3$ in the class group of $X$.
This condition is very restrictive, for example, it implies that $X$ has at most $2$ irreducible components \cite[7.1]{H04}.
(On the other hand, already for $d=6$ there are surfaces $X$ with $18$ components.)

\begin{remark}
It is instructive to compare with the one dimensional case.
We consider the moduli space $\cM$ of pairs $(X,D)$ such that $X \simeq \bP^1$ and $D$ is a divisor of degree $d \ge 3$. Let us order the points $D=P_1+\cdots+P_d$ for simplicity. We construct a compactification $\ocM$ as a moduli space of pairs $(X,D)$ such that
$(X,(\frac{2}{d}+\epsilon)D)$ is slc and $K_X+(2/d+\epsilon)D$ is ample for $0 < \epsilon \ll 1$. The first condition means simply that $X$ is a nodal curve, the points $P_i$ are smooth points of $X$, and the divisor $(\frac{2}{d}+\epsilon)D$ has coefficients $\le 1$ .
Then $dK_X+ 2D \sim 0$, and $X$ is either a copy of $\bP^1$, or two copies of $\bP^1$ meeting transversely in a single point.
Moreover, in the second case $d$ is even and there are $d/2$ points on each component. We deduce that there is a birational morphism $\ocM \rightarrow \ocM^{\GIT}$ to the symmetric GIT quotient
$\ocM^{\GIT}=(\bP^1)^d/\!/\SL(2)$ which is an isomorphism if $d$ is odd and is a resolution of the singularities of $\ocM^{\GIT}$ corresponding to the strictly semistable points if $d$ is even.
\end{remark}

We describe the degenerate surfaces $X$ in case $d$ is not divisible by $3$. There are two types $A$ and $B$.
Surfaces of type $A$ are normal and have quotient singularities of Wahl type.
Surfaces of type $B$ have two irreducible components meeting along a smooth rational curve.
We have a complete description of the surfaces of type $A$:
\begin{theorem}\label{typeA}\cite[1.2]{HP10}
Let $X$ be a normal surface with quotient singularities which admits a smoothing to the projective plane.
Then $X$ is obtained from a weighted projective plane $\bP(a^2,b^2,c^2)$
by a $\bQ$-Gorenstein deformation that smoothes some subset of its singularities, where $(a,b,c)$ is a solution of the Markov equation
$$a^2+b^2+c^2=3abc.$$
\end{theorem}
The solutions of the Markov equation are easily described: $(1,1,1)$ is a solution, and all solutions are obtained from $(1,1,1)$ by a sequence of \emph{mutations} of the form
\begin{equation}\label{mutation}
(a,b,c) \mapsto (a,b,c'=3ab-c).
\end{equation}
We can define a graph $G$ with vertices labelled by solutions of the Markov equation and edges corresponding to pairs of solutions related by a single mutation. Then $G$ is an infinite tree such that every vertex has degree $3$.

The surface $\bP(a^2,b^2,c^2)$ has cyclic quotient singularities of types $\frac{1}{a^2}(b^2,c^2)$, $\frac{1}{b^2}(c^2,a^2)$, $\frac{1}{c^2}(a^2,b^2)$.
Using the Markov equation one sees that these are Wahl singularities (note that $a,b,c$ are coprime and not divisible by $3$ by the inductive description of the solutions of the Markov equation above). Moreover there are no equisingular deformations and no local-to-global obstructions because $H^1(\cT_{\bP})=H^2(\cT_{\bP})=0$. So there is one $\bQ$-Gorenstein deformation parameter for $\bP$ associated to each singularity.
A surface $X$ obtained as a $\bQ$-Gorenstein deformation of $\bP$ is uniquely determined by the subset of singularities that is smoothed.

If $X$ is a surface of type $B$ then $X$ is the union of two normal surfaces $X_1$ and $X_2$ meeting along a smooth rational curve $C$.
The surface $X$ has orbifold normal crossing singularities $(xy=0) \subset \frac{1}{n}(1,-1,a)$ along the double curve $C$ and Wahl singularities away from $C$. The Picard numbers of the components of $X$ are given by either $\rho(X_1)=\rho(X_2)=1$ or $\{\rho(X_1),\rho(X_2)\}=\{1,2\}$.

\begin{example}\label{basicex}
Consider the two parameter family
$$\cX=(X_0X_2=sX_1^2+tY) \subset \bP(1,1,1,2) \times \bC^2_{s,t}$$
The special fiber $X:=\cX_0$ is a surface of type $B$ obtained by glueing two copies of the weighted projective plane $\bP(1,1,2)$ (the quadric cone) along a ruling. It has an orbifold normal crossing singularity of type $(xy=0) \subset \frac{1}{2}(1,1,1)$.
The fibers $\cX_{s,t}$ for $s \neq 0$, $t=0$ are isomorphic to $\bP(1,1,4)$, with the embedding in $\bP(1,1,1,2)$ being the $2$-uple embedding
$$\bP(1,1,4) \rightarrow (X_0X_2=X_1^2) \subset \bP(1,1,1,2)$$
$$(U_0,U_1,V) \mapsto (X_0,X_1,X_2,Y)=(U_0^2,U_0U_1,U_1^2,V).$$
The surface $\bP(1,1,4)$ has a Wahl singularity of type $\frac{1}{4}(1,1)$.
The fibers $\cX_{s,t}$ for $t \neq 0$ are isomorphic to $\bP^2$.
The deformation $\cX/(0 \in \bC^2_{s,t})$ is the versal $\bQ$-Gorenstein deformation of $X$.
\end{example}

We remark that in general it is difficult to describe the versal $\bQ$-Gorenstein deformation of $\bP:=\bP(a^2,b^2,c^2)$ explicitly, because for example the embedding in weighted projective space defined by $-\frac{1}{3}K_{\bP}=\cO_{\bP}(abc)$ (corresponding to $\cO_{\bP^2}(1)$ on the general fiber $\bP^2$) has high codimension.

\begin{example}
Here we describe a two parameter family of surfaces which ``connects'' the weighted projective planes $\bP:=\bP(a^2,b^2,c^2)$, $\bP':=\bP(a^2,b^2,{c'}^2)$ associated to two solutions of the Markov equation related by a single mutation (\ref{mutation}).
The family is given by
$$\cX=(XY=sZ^{c'}+tT^c) \subset \bP(a^2,b^2,c,c') \times \bC^2_{s,t}.$$
The special fibre $X:=\cX_0$ is the union of two weighted projective planes $\bP(a^2,c,c')$, $\bP(b^2,c,c')$ glued along the coordinate lines of degree $a^2$ and $b^2$. It has two Wahl singularities of indices $a$ and $b$ and two orbifold normal crossing singularities of indices $c$ and $c'$.
The fibers $\cX_{s,t}$ for $s = 0$, $t \neq 0$ are isomorphic to $\bP(a^2,b^2,c^2)$, with the embedding being the $c$-uple embedding
$$\bP(a^2,b^2,c^2) \rightarrow (XY=T^c) \subset \bP(a^2,b^2,c,c')$$
$$(U , V , W) \mapsto (X , Y , Z, T) = (U^{c},V^{c},W,UV).$$
Similarly, the fibers $\cX_{s,t}$ for $s \neq 0$, $t=0$ are isomorphic to $\bP(a^2,b^2,{c'}^2)$.
The fibers $\cX_{s,t}$ for $s \neq 0$, $t \neq 0$ are obtained from $\bP$ or $\bP'$ by smoothing the singularity of index $c$ or $c'$ respectively.
Note that Example~\ref{basicex} is the special case $a=b=c=1$.
\end{example}

\begin{example}
Here we describe surfaces $X$ of type $B$ such that $\rho(X_1)=2$, $\rho(X_2)=1$.
We begin with the trivial family $\cY=\bP^2 \times \bC^1_t$.
Let $m,n$ be positive integers such that $(m,n)=1$.
Let $\pi \colon \cX \rightarrow \cY$ be the blowup of the point
$P=((1 \colon 0 \colon 0),0)$ with weights $(m,n,1)$ with respect to the coordinates $x_1=X_1/X_0,x_2=X_2/X_0,t$.
The special fiber $X:=\cX_0$ is reduced and has two components $X_1$ and $X_2$ given by the strict transform of $Y:=\cY_0$ and the exceptional divisor $E=\bP(m,n,1)$ of $\pi$ respectively.
The restriction $p \colon X_1 \rightarrow Y=\bP^2$ is the weighted blowup of the point $(1\colon 0 \colon 0)$ in $\bP^2$ with weights $(m,n)$ with respect to the coordinates $x_1,x_2$. The components $X_1$ and $X_2$ are glued along the smooth rational curve $C$ given by the exceptional divisor of $p$ on $X_1$ and the coordinate line of degree $1$ on $X_2=\bP(m,n,1)$. The surface $X$ has two orbifold normal crossing singularities of indices $m$ and $n$.
A toric calculation shows that the $\bQ$-Cartier divisor $-K_X$ is ample iff $\frac{1}{2}m < n < 2m$.
In this case the surface $X$ occurs in $\ocP_d$ for $d$ divisible by $3mn$.
The canonical divisor $-K_X$ is divisible by $3$ in the class group of $X$ iff $m+n$ is divisible by $3$.
In this case $X$ occurs in $\ocP_d$ for $d$ divisible by $mn$.

\end{example}

\subsection{Boundary divisors}
Let $\cP'_d \subset \ocP_d$ denote the locus of pairs $(X,D)$ such that the surface $X$ is isomorphic to $\bP^2$, and let $\partial \ocP_d := \ocP_d \setminus \cP'_d$. Suppose $d$ is not a multiple of $3$. Then $\partial\ocP_d \subset \ocP_d$ is a normal crossing divisor with irreducible components corresponding to surfaces of the following types.
\begin{enumerate}
\item Surfaces of type $A$ with a unique Wahl singularity.
\item Surfaces $X=X_1 \cup X_2$ of type $B$ such that $\{\rho(X_1),\rho(X_2)\}=\{1,2\}$ and $X$ is smooth away from $C:=X_1 \cap X_2$.
\end{enumerate}
This follows from an analysis of the versal $\bQ$-Gorenstein deformations of surfaces of type $A$ and $B$, together with the following fact: for
$[(X,D)] \in \ocP_d$, the forgetful map $\Def^{\QG}(X,D) \rightarrow \Def^{\QG}(X)$ from deformations of the pair to deformations of the surface is smooth. See \cite[3.12, 8.2, 9.1]{H04}.
(Note: The surfaces $X$ which occur in a given degree $d$ are those for which there exists a divisor $D \in |-(d/3)K_X|$ such that the pair $(X,(3/d+\epsilon))$ is slc for $0 < \epsilon \ll 1$. A necessary condition is that the index of each singularity $P \in X$ is at most $d$ \cite[4.5]{H04}.)

\begin{remark}
Our stability condition for a pair $(X,D)$ such that $X \simeq \bP^2$ is a natural strengthening of GIT stability for the plane curve $D$ to a local analytic condition \cite[\S10]{H04}. In particular, the locus $\cP'_d$ is contained in the moduli stack of GIT stable plane curves of degree $d$.
\end{remark}

\section{Exceptional vector bundles associated to degenerations of surfaces}\label{bundles}

Let $Y$ be a projective surface. A vector bundle $F$ on $Y$ is \emph{exceptional} if $\End F = \bC$ and $H^1(\cEnd F)=H^2(\cEnd F)=0$.
In particular, an exceptional vector bundle $F$ is indecomposable, rigid (no infinitesimal deformations), and unobstructed in families.
So, if $\cY/(0 \in S)$ is a deformation of $Y$ over a germ $(0 \in S)$, then there exists a unique vector bundle $\cF$ on $\cY$ such that $\cF|_Y=F$.

\begin{theorem} \label{mainthm} \cite{H11}
Let $X$ be a projective normal surface with a unique singularity $P \in X$ of type $\frac{1}{n^2}(1,na-1)$.
Let $\cX/(0 \in T)$ be a one parameter $\bQ$-Gorenstein deformation of $X$ such that the general fiber is smooth.
Let $Y$ denote a general fiber of $\cX/T$.
\begin{enumerate}
\item Assume that $H_1(Y,\bZ)$ is finite of order coprime to $n$.
Then the specialization map
$$\spz \colon H_2(Y,\bZ) \rightarrow H_2(X,\bZ)$$
is injective with cokernel isomorphic to $\bZ/n\bZ$.
\item Assume in addition that $H^{2,0}(Y)=0$. Then, after a base change $T' \rightarrow T$ of degree $a$, there exists a reflexive sheaf $\cE$ on $\cX':=\cX \times_T T'$ satisfying the following properties.
\begin{enumerate}
\item $F:=\cE|_Y$ is an exceptional vector bundle of rank $n$ on $Y$.
\item $E:=\cE|_X$ is a torsion-free sheaf on $X$ such that its reflexive hull $E^{\vee\vee}$ is isomorphic to the direct sum of $n$ copies of a reflexive rank $1$ sheaf $A$, and the quotient $E^{\vee\vee}/E$ is a torsion sheaf supported at $P \in X$.
\end{enumerate}
If $\sH$ is a line bundle on $\cX/T$ which is ample on fibers, then the vector bundle $F$ is slope stable with respect to $H:=\sH|_Y$.
The Chern classes of $F$ are given by
$$c_1(F)=n c_1(A) \in H_2(Y,\bZ) \subset H_2(X,\bZ),$$
$$c_2(F)=\frac{n-1}{2n}(c_1(F)^2+n+1).$$
Moreover
$$c_1(F) \cdot K_Y = \pm a \mod n,$$
and
$$H_2(X,\bZ)=H_2(Y,\bZ)+\bZ\cdot(c_1(F)/n).$$
\end{enumerate}
\end{theorem}

\begin{remark}
Note that the bundles obtained from $F$ by dualizing or tensoring by a line bundle arise in the same way. (Indeed, if $L$ is a line bundle on $Y$ then $L$ extends to a reflexive rank one sheaf $\sL$ on $\cX$. Now $\cE^{\vee}$ and $(\cE \otimes \sL)^{\vee\vee}$ satisfy the properties \ref{mainthm}(2a,b) and restrict to $F^{\vee}$ and $F \otimes L$ on $Y$.)
\end{remark}

We sketch the proof of Theorem~\ref{mainthm} in the case $n=2$. In this case $P \in X$ is a singularity of type $\frac{1}{4}(1,1)$.
Let us assume for simplicity that the $\bQ$-Gorenstein deformation $(P \in \cX)/(0 \in T)$ is versal, that is, isomorphic to
$$(xy=z^2+t) \subset \left( \bC^3_{x,y,z} / \textstyle{\frac{1}{2}}(1,1,1) \right) \times \bC^1_t.$$
(In general, we obtain our construction by pullback from the versal case.)

The germ $(P \in \cX)$ is a cyclic quotient singularity of type $\frac{1}{2}(1,1,1)$.
Let $\pi \colon \tilde{\cX} \rightarrow \cX$ denote the blowup of $P \in \cX$.
Then the exceptional locus $W$ of $\pi$ is a copy of the projective plane with normal bundle $\cO(-2)$.
The total space $\tilde{\cX}$ is smooth and the special fiber $\tilde{X}:=\tilde{\cX}_0$ is a normal crossing divisor
with two smooth components given by the strict transform $X'$ of $X \subset \cX$ and the exceptional divisor $W$.
The restriction $p \colon X' \rightarrow X$ of $\pi$ is the minimal resolution of $(P \in X)$.
The surfaces $X'$ and $W$ meet along the smooth rational curve $C$  which is the exceptional curve of $p \colon X' \rightarrow X$ and is embedded as a conic in $W \simeq \bP^2$.

Let $B$ denote the intersection of $X$ with a small ball with center at the singularity $P \in X$ in some embedding.
Let $L$ denote the boundary of $B$ (the link of the singularity) and write $X^o := X \setminus B$.
Our assumption on $H_1(Y,\bZ)$ implies that the map
$$H^2(X^o,\bZ) \rightarrow H^2(L,\bZ)$$
is surjective, by a local analysis near the singularity (see \S\ref{wahl}) and a Mayer--Vietoris argument.
(See \cite{H11}, proof of Theorem~1.1, for more details.)
It follows that the restriction map
\begin{equation}\label{surjection}
H^2(X',\bZ) \rightarrow H^2(C,\bZ)
\end{equation}
is surjective.

We have
$$H^i(\cO_{X'})=H^i(\cO_X)=H^i(\cO_Y)$$
because the singularity $P \in X$ is a quotient singularity, and by assumption $H^i(\cO_Y)=0$ for $i > 0$.
Thus $H^i(\cO_{X'})=0$ for $i>0$ and the map
$$c_1 \colon \Pic X' \rightarrow H^2(X',\bZ)$$
is an isomorphism.
So by the surjectivity of (\ref{surjection}) above, there is a line bundle $A'$ on $X'$ such that the restriction of $A'$ to $C$ has degree $1$.

Let $G$ be the vector bundle $\cT_{\bP^2}(-1)$ on $W=\bP^2$. Then $G$ is an exceptional vector bundle on $W$, and $G|_C \simeq \cO_C(1)^{\oplus 2}$, where $\cO_C(1)$ denotes the line bundle of degree $1$ on the smooth rational curve $C$.
Since $\tilde{X}=X'\cup W$ is a normal crossing surface with double curve $C= X' \cap W$, we have the exact sequence of sheaves on $\tilde{X}$
\begin{equation}\label{devissage}
0 \rightarrow \cO_{\tilde{X}} \rightarrow \cO_{X'} \oplus \cO_W \rightarrow \cO_C \rightarrow 0.
\end{equation}
Thus we can define a vector bundle $\tilde{E}$ over $\tilde{X}$ by glueing $A'^{\oplus 2}$ over $X'$ and $G$ over $W$ along $\cO_C(1)^{\oplus 2}$ over $C$.

We show that $\tilde{E}$ is exceptional. Tensoring the exact sequence (\ref{devissage}) with $\cEnd \tilde{E}$ we obtain the exact sequence
$$0 \rightarrow \cEnd \tilde{E} \rightarrow \cO_{X'}^{2 \times 2} \oplus \cEnd G \rightarrow \cO_C^{2 \times 2} \rightarrow 0.$$
Now $H^1(\cO_C)=0$ because $C$ is a smooth rational curve and $H^i(\cO_{X'})=0$ for $i>0$ as noted above, so we obtain $H^i(\cEnd \tilde{E}) = H^i(\cEnd G)$. Thus $\tilde{E}$ is exceptional because $G$ is exceptional.

Let $\tilde{\cE}$ denote the (unique) vector bundle over $\tilde{\cX}$ obtained by deforming the exceptional bundle $\tilde{E}$.
Then the restriction $F$ of $\tilde{\cE}$ to a general fiber $Y$ of $\cX/T$ is exceptional by upper semicontinuity of cohomology.
Let $\cE:=(\pi_*\cE)^{\vee\vee}$ be the reflexive hull of the pushforward of $\cE$ to $\cX$.
Similarly let $A:=(p_*A')^{\vee\vee}$ be the reflexive hull of the pushforward of $A'$ to $X$.
Then $E:=\cE|_X$ is torsion-free because $\cE$ is reflexive, and by construction $E|_{X \setminus \{P\}}=A^{\oplus 2}|_{X \setminus \{P\}}$, so $E^{\vee\vee}=A^{\oplus 2}$ and $E^{\vee\vee}/E$ is supported at $P$.

\section{Boundary divisors of the moduli space of stable surfaces}\label{conjecture}

Here we formulate a precise conjecture relating the components of the boundary of the moduli space of stable surfaces and exceptional vector bundles.

Let $Y$ be a smooth surface such that $K_Y$ is ample, $H^2(\cT_Y)=0$, $H^{2,0}(Y)=0$, and $\pi_1(Y)=0$.
These conditions imply that $1 \le K_Y^2 \le 5$ and $Y$ is homeomorphic (but not diffeomorphic) to the blowup of $\bP^2$ in $9-K_Y^2$ points.
Examples of such surfaces $Y$ are known for $K_Y^2=1,2,3,4$ \cite{B85}, \cite{LP07}, \cite{PPS09a}, \cite{PPS09b}.

Let $\ocM$ denote the irreducible component of the moduli stack of stable surfaces containing $[Y]$.
Let $\cM \subset \ocM$ denote the locus of surfaces with Du Val singularities and $\partial \ocM := \ocM \setminus \cM$.
Note that $\cM$ is smooth of the expected dimension $10-2K_Y^2$ near $[Y]$ by our assumption $H^2(\cT_Y)=0$.
If $\ocM$ is not normal, we replace it by its normalization.

Let $\Gamma \subset \Aut(H^2(Y))$ denote the monodromy group.
Then $\Gamma$ preserves $K_Y$.
The lattice $K_Y^{\perp}$ is negative definite because $H^{2,0}=0$ and $K_Y$ is ample.
Thus $\Gamma$ is a finite group.

For $F$ a vector bundle on $Y$ we define its \emph{slope vector} $v(F):=c_1(F)/\rk(F) \in H^2(Y,\bQ)$.
Then $v(F \otimes L) = v(F) + c_1(L)$ and $v(F^{\vee})=-v(F)$. Note that $c_1 \colon \Pic(Y) \rightarrow H^2(Y,\bZ)$ is an isomorphism because $H^1(\cO_Y)=H^2(\cO_Y)=0$ by assumption.

Consider exceptional vector bundles $F$ of rank greater than $1$ on $Y$ that are slope stable with respect to $K_Y$ and of degree $c_1(F) \cdot K_Y$ coprime to $\rk(F)$.
Let $S$ denote the set of associated slope vectors $v(F)=c_1(F)/\rk(F) \in H^2(Y,\bQ)$ modulo translation by $H^2(Y,\bZ)$, multiplication by $\pm 1$, and the action of the monodromy group $\Gamma$.

Let $T$ denote the set of codimension one components $D$ of the boundary $\partial \ocM$ such that the general fiber over $D$ is a normal surface $X$ with a unique singularity $P \in X$ of Wahl type.

Theorem~\ref{mainthm} defines a map of sets

$$\Phi \colon T \rightarrow S, \quad D \mapsto v(F) .$$
We hope that the map $\Phi$ is ``close" to a bijective correspondence.
In particular, we formulate the following conjecture.

\begin{conjecture}
The set $S$ is finite. Equivalently, the ranks of exceptional bundles $F$ on $Y$ which are slope stable with respect to $K_Y$ (and of degree coprime to the rank) are bounded.
\end{conjecture}

One can show that the analogue of the map $\Phi$ for $Y=\bP^2$ is bijective using the classification of degenerations of $\bP^2$ given by Theorem~\ref{typeA} and the classification of exceptional vector bundles on $\bP^2$ \cite{DLP85},\cite{R89}. More precisely:

\begin{theorem}\cite[\S6]{H11}
There is a bijective correspondence between isomorphism types of degenerations $X$ of $\bP^2$ with a unique Wahl singularity and isomorphism types of exceptional bundles $F$ of rank greater than $1$ on $\bP^2$ modulo dualizing and tensoring by line bundles, given by the construction of Theorem~\ref{mainthm}.
\end{theorem}

\section{Relation with Donaldson theory}\label{donaldson}

We note that the proposal of \S\ref{conjecture} is closely related to the Donaldson theory of invariants of smooth $4$-manifolds \cite{DK90}, \cite{K05}, \cite{M09}.

Donaldson showed that the classification of simply connected smooth $4$-manifolds up to diffeomorphism is much richer than the classification up to homeomorphism obtained by Freedman \cite{F82}. This was achieved by defining new invariants of smooth $4$-manifolds as follows. Let $X$ be a simply connected compact oriented smooth $4$-manifold.
Let $b_2^+$ denote the dimension of a maximal positive definite subspace of $H^2(X,\bR)$; we assume $b_2^+>0$. (If $X$ is a K\"ahler manifold, then $b_2^+=2\dim H^{2,0}+1$ by Hodge theory.)
Fix a Riemannian metric $g$ on $X$ and consider a smooth complex vector bundle $E$ of rank $n$ on $X$ together with a Hermitian metric on $E$.
(Note that the isomorphism type of $E$ is determined by its Chern classes because $X$ has real dimension $4$ \cite{W52}.)
One considers the moduli space $M$ of \emph{instantons}: anti-self-dual connections on the principal $\PU(n)$-bundle associated to $E$ (where $\PU(n)=\U(n)/\U(1)$ denotes the projective unitary group). One defines numerical invariants of $(X,g)$ as intersection numbers of certain natural cohomology classes on $M$.
If $b_2^+>1$ then the Donaldson invariants are independent of the choice of the metric $g$, that is, they are invariants of the smooth manifold $X$. If $b_2^+=1$ then they depend on the metric $g$ through the associated line of self-dual harmonic $2$-forms in $H^2(X,\bR)$ via a chamber structure and wall crossing formulas.

If now $X$ is a complex surface and $g$ is a K\"ahler metric, then the moduli space of irreducible anti-self-dual connections is identified with the moduli space of holomorphic structures on $E$ which are slope stable with respect to the class of the K\"ahler form. Thus, in the case of an algebraic surface, Donaldson invariants can be computed by intersection theory on the moduli space of stable vector bundles.
Moreover, in case $b_2^+=1$, the line of self dual harmonic $2$-forms is spanned by the K\"ahler form, and the wall crossing formulas are determined by the variation of the moduli of stable bundles with the polarization.

One of the first applications of Donaldson theory was the analysis of connected sum decompositions of smooth $4$-manifolds.
Let $X=X_1 \# X_2$ be a connected sum decomposition of a $4$-manifold $X$.
Thus $X=X_1^o \cup X_2^o$ where $X_i^o$ is the complement of a ball about a point in $X_i$, and $X$ is obtained by identifying the boundaries of $X_1^o$ and $X_2^o$.
We choose a Riemannian metric $g$ on $X$ such that the ``neck'' connecting $X^o_1$ and $X^o_2$ is isometric to $S^3 \times (-R,R)$ for $R \gg 0$.
In this situation one proves that instantons on $X$ are obtained by gluing instantons on $X^o_1$ and $X^o_2$ which decay exponentially to a flat connection at the boundary. In this way one obtains gluing formulas for the Donaldson invariants of $X$ in terms of relative Donaldson invariants of $X_1^o$ and $X_2^o$.

The situation we consider above is similar. Recall that we begin with a smooth surface $Y$ and consider a degeneration to a normal surface $X$ with a unique Wahl singularity $P \in X$. We have a \emph{generalized connected sum decomposition} $Y = X^o \#_L M$ of the $4$-manifold $Y$. That is, writing $X^o$ for the complement of a small ball about the singular point $P \in X$, $L$ for the link of the singularity, and $M$ for the Milnor fiber of the smoothing, the $4$-manifold $Y$ is obtained by identifying the boundaries of $X^o$ and $M$, which are copies of the lens space $L$.
Such generalized connected sum decompositions may be studied using Donaldson theory in the same way as above \cite{D02}.
In particular, the rank $2$ Donaldson invariants for the decomposition corresponding to a Wahl singularity of type $\frac{1}{n^2}(1,n-1)$ were analyzed in \cite{FS97}. (We note that the surgery of smooth 4-manifolds given by passing from the minimal resolution of a Wahl singularity to its smoothing is called a \emph{rational blowdown} in differential and symplectic geometry.)

We expect that our construction admits an interpretation in these terms.
More precisely, we expect that for an appropriate choice of metric $g$, the exceptional bundle we construct admits a unique anti-self-dual connection which is projectively trivial over $X^o$ (and nontrivial over $M$).

In the Donaldson theory it is common to restrict attention to bundles of rank $2$. In fact it is conjectured that the higher rank invariants are determined by the rank $2$ invariants \cite{MM98}.
However, as we have seen, for our application to the study of the boundary of the moduli space of complex surfaces, it is crucial to consider vector bundles of arbitrary rank $n$, because they are directly connected to the geometry of the degeneration. On the other hand, we restrict our attention to rigid vector bundles, that is, we only consider the case in which the moduli space of stable vector bundles is zero dimensional.

\section{Other examples of boundary divisors}

We describe some types of codimension $1$ boundary components of the moduli space of stable surfaces different from the Wahl type described in \S4.

Let $X$ be a stable surface. If $X$ satisfies the conditions set out in one of the examples below, then equisingular deformations of $X$ determine a codimension $1$ component of the boundary of the moduli space of stable surfaces.

\begin{example}
The surface $X$ has a unique singularity $P \in X$ of one of the following types, and there are no local-to-global obstructions.
\begin{enumerate}
\item The minimal resolution $\pi \colon \tilde{X} \rightarrow X$ has exceptional locus a union of $4$ smooth rational curves $E_1,E_2,E_3,$ and $F$.
The $E_i$ are disjoint, each meets the curve $F$ transversely in a single smooth point, and
$$(-E_1^2,-E_2^2,-E_3^2;-F^2)=(3,3,3;4),(2,4,4;3),\mbox{ or }(2,3,6;2).$$
See \cite[\S8]{W10}.  (Note: The index one cover $p \colon (Q \in Z) \rightarrow (P \in X)$ is the quotient of the cone over an elliptic curve of degree $9$, $8$, or $6$ by a cyclic group of order $3$, $4$, or $6$, respectively \cite[9.6(3)]{K88}.)
\item A cone over an elliptic curve of degree $9$ \cite[4.5]{S98}.
\item A smoothable cusp singularity such that, writing $E=E_1+\cdots+E_r$ for the exceptional locus of the minimal resolution, $-E^2=9+r$ \cite[5.6]{W81}.
(Note: The smoothability criterion for cusp singularities conjectured by Looijenga \cite[III.2.11]{L81} is proven in \cite{GHK11}. The smoothable cusps with $r \le 3$ are listed in \cite[\S4]{FM83}.)
\end{enumerate}
\end{example}
\begin{example}
The surface $X$ is the union of two components $X_1$ and $X_2$ meeting transversely along a smooth curve $C$, $X$ has orbifold normal crossing singularities $(xy=0) \subset \frac{1}{n}(1,-1,a)$ along $C$ and is smooth elsewhere. Then $\cT^1_{\QG,X}$ is the line bundle on $C$ given by (\ref{QGinfnorbdle}). We require that $\dim H^0(\cT^1_{\QG,X})=1$, $H^1(\cT^1_{\QG,X})=0$, and $H^2(\cT_X)=0$. For example if $C$ has genus zero we require $\cT^1_{\QG,X} \simeq \cO_C$. (See \cite{T09} for more results on $\bQ$-Gorenstein deformations of non-normal surfaces.)
\end{example}


\bibliographystyle{amsplain}

\end{document}